\documentclass[11pt,a4paper]{article}

\usepackage[utf8]{inputenc}
\usepackage[T1]{fontenc}
\usepackage[english]{babel}
\usepackage{lmodern}
\usepackage{microtype}
\usepackage{geometry}
\usepackage{setspace}
\usepackage[hidelinks]{hyperref}
\usepackage{enumitem}
\usepackage{booktabs}
\usepackage{csquotes}
\usepackage[round,authoryear]{natbib}

\title{\textbf{Louis Bachelier's \emph{Th\'eorie de la sp\'eculation}:\\
Circulation and Reception}}

\author{Hector Giacomini\\
\small Institut Denis Poisson, Universit\'e de Tours--CNRS, Tours, France\\
\small\texttt{giacominihector@gmail.com}}

\date{}
\usepackage{amsfonts}

\begin{document}
\maketitle

\begin{abstract}
Louis Bachelier's \emph{Théorie de la spéculation} (1900)
introduced the first probabilistic model of Brownian motion
in a study of financial fluctuations, five years before
Einstein's physical theory of 1905. By relating the evolution
of probability distributions to the heat equation, Bachelier
anticipated central elements of Brownian-motion theory,
which would become a major field of research in mathematics
and physics. His treatment of speculation and option valuation
also laid foundations for modern mathematical finance.
The subsequent importance of these contributions raises
a historical question: how widely did the thesis circulate
before its broader significance was recognized?

Drawing on correspondence with libraries, surviving copies,
accession registers, provenance marks, historical catalogues,
and printed acquisition lists, this article reconstructs
the international circulation and reception of the thesis
in its two published forms: the separate Gauthier-Villars
dissertation and its publication in the \emph{Annales
scientifiques de l'École normale supérieure}. It identifies
the thirty-eight foreign academic institutions listed in the official
dissertation exchange network for 1900 and documents its
operation through surviving copies and records of receipt.
University exchanges, commercial purchases, journal
distribution, and personal presentations provided
complementary channels of dissemination. Extensive holdings
in American university libraries reveal its geographical
reach, with dated evidence establishing receipt at Cornell
in April 1900. Publisher promotion and repeated authorial
donations further document efforts to reach readers across
disciplinary and national boundaries.

Contemporary reviews and subsequent uses establish selective
engagement in mathematics, actuarial science, and economics,
including transmission without attribution. An examination
of probability composition and diffusion clarifies the
thesis's place in the history of stochastic processes
and its mathematical relationship to later treatments
of physical Brownian motion. Surviving copies also connect
early dissemination with renewed interest: Columbia's
journal volume, accessioned in 1901, records Mandelbrot's
borrowing in 1961.

Distinguishing material availability from documented reading
and use, the study shows that delayed recognition coexisted
with early international circulation and an identifiable,
though dispersed, readership. It thereby extends previous
reassessments of Bachelier's reception and examines how
ideas developed for financial speculation came to be
appreciated within different scientific traditions.
\end{abstract}

\section{Bachelier's thesis: academic background, subject, and reception}

Bachelier's route into mathematics differed from that
of much of the French scientific elite. Born in Le Havre in 1870, he interrupted his studies
after both parents died in 1889 and had to work.
He may have taken over the family business. Military service further delayed
his education, and he began studying at the Sorbonne
in 1892, aged twenty-two. He had attended neither
the \'Ecole normale sup\'erieure nor the \'Ecole polytechnique,
which provided structured training and established routes
into teaching, administration, and the state's technical
services. His university studies thus fell outside
the institutional framework that secured early professional
prospects for many future members of the scientific elite.
During these years, he appears to have worked in connection
with the Paris Bourse; his dissertation demonstrates
detailed familiarity with contemporary trading practices
\citep[p.~4]{Taqqu2001}.

His training nevertheless brought him into direct contact
with leading mathematicians, including Paul Appell, Joseph
Boussinesq, and Henri Poincar\'e \citep{Courtault2000}.
As Davis and Etheridge emphasize, Boussinesq introduced
him to heat theory and Poincar\'e to probability,
providing resources he would combine in his analysis
of financial fluctuations
\citep[pp.~80--81]{DavisEtheridge2006}.

Bachelier defended his thesis \emph{Th\'eorie de la
sp\'eculation} at the Faculty of Sciences of Paris on
29 March 1900, before a commission comprising Appell
as president, with Boussinesq and Poincar\'e as examiners
\citep{Courtault2000}.
It was published separately by Gauthier-Villars and,
in the same year, in volume 17 of the third series of
the \emph{Annales scientifiques de l'\'Ecole normale
sup\'erieure}, pp.~21--86
\citep{Bachelier1900,Bachelier1900Separate}.
The separate printing also served an institutional purpose:
the doctoral procedure required candidates to provide
the University with a substantial number of copies
for examining bodies and institutional distribution,
including exchanges with French and foreign university
libraries.

The thesis applies the calculus of probabilities to
speculation on the Paris Bourse, particularly transactions
in French government bonds known as \emph{rentes}. Before developing his mathematical analysis, Bachelier introduces
the principal stock-exchange operations and explains the terminology
used in their description \citep[p.~6]{Bachelier1900Separate}.

Participants formed judgments about future prices
by interpreting market conditions. In his introduction,
Bachelier emphasized the many influences affecting prices
and the conflicting expectations of buyers and sellers.
He sought to express mathematically the probabilities
implicit in the market's assessment of possible price
variations.
Earlier probabilistic approaches to speculation existed,
notably that of Jules Regnault in 1863
\citep[pp.~13--14]{Taqqu2001}.
Bachelier's distinctive contribution was a systematic
mathematical treatment of price fluctuations in continuous
time and its application to financial contracts.

Its historical importance extends far beyond this setting.
Bachelier described price variations through probability
laws evolving continuously with time. He connected
the propagation of probability with Fourier's heat equation,
obtained Gaussian distributions whose characteristic scale
grows as the square root of time, and calculated
the probability of reaching a specified price level
within a given interval. The mathematical significance
of this construction will be examined below.

Poincar\'e's report combines caution about the unusual
subject with appreciation of the work's mathematical
originality. He noted that speculation lay outside
the subjects ordinarily treated by doctoral candidates,
while recognizing Bachelier's efforts to delimit
the legitimate application of probability to financial
transactions. Addressing the thesis's central analytical
argument, he considered the analogy between the propagation
of probability and heat legitimate: Fourier's reasoning
could be applied almost unchanged to this different problem
\citep{Poincare1900}.
Although he regretted that Bachelier had not developed
the connection further, he judged it sufficient to justify
the Gaussian law and indicate its limits of applicability
\citep[p.~78]{DavisEtheridge2006}.

Poincar\'e also singled out Bachelier's solution
of a difficult first-passage problem: determining
the probability that a price would reach a specified
level before a given date. He praised the argument
as short, simple, elegant, and original \citep{Poincare1900}.

Despite this favourable assessment, the thesis received
the mention \emph{honorable}, rather than the higher
distinction \emph{tr\`es honorable}.
This outcome requires careful interpretation.
As \citet{Courtault2000} emphasized, Poincar\'e's attentive
and favourable report makes it difficult to describe
the thesis as simply dismissed or misunderstood
by its examiners. It entered circulation as a recognized
mathematical dissertation, examined by leading French
mathematicians and published in one of the country's
principal mathematical journals.

The doctorate did not secure Bachelier an academic position.
Between 1909 and 1914, he gave a \emph{cours libre}
at the Faculty of Sciences of Paris on probability,
financial operations, and analogies with physics,
without receiving remuneration for this teaching.

In 1914, the University Council supported efforts
to make his appointment permanent, but his mobilization
interrupted this prospect. He served in the French army
throughout the war and never obtained a permanent
professorship in Paris. After the war, he held temporary
appointments as a \emph{chargé de cours}, replacing
professors on leave or filling vacancies pending
a permanent appointment, at Besan\c{c}on (1919--1922),
Dijon (1922--1925), and Rennes (1925--1927)
\citep[p.~17]{Taqqu2001}.
Only in 1927, aged fifty-seven, did he become a full
professor at Besan\c{c}on, remaining until his retirement
in 1937 \citep{Courtault2000,BenElMechaiekhDimand2018}.

This prolonged professional insecurity forms part
of the context in which his work circulated, but does
not establish that the thesis remained unavailable
or unread. Academic recognition, publication distribution,
and the use of mathematical results followed different
trajectories. Understanding its early history therefore
requires examining how the thesis reached potential readers
and what evidence survives of their engagement with it.

Previous historical studies have challenged the image
of Bachelier as an entirely forgotten precursor,
highlighting contemporary assessments, early readers,
and subsequent uses of his ideas
\citep{Courtault2000,Taqqu2001,Jovanovic2002,
DavisEtheridge2006,Jovanovic2012,BenElMechaiekhDimand2018}.
The present study contributes to this reassessment
by reconstructing the material and institutional
circulation of the thesis in its two published forms:
the separate Gauthier-Villars edition and its appearance
in the \emph{Annales scientifiques de l'\'Ecole normale
sup\'erieure}. It combines surviving copies, provenance
marks, accession records, historical catalogues,
publishers' announcements, and printed donation lists
to establish where and when the text's presence
can be documented and how it reached different institutions.

These sources reveal complementary routes of dissemination:
university dissertation exchanges, commercial purchases,
personal presentations, and the journal's international
distribution. They also document Bachelier's efforts
to reach readers outside France, including members
of the actuarial profession. This international circulation
began in the year of publication and continued through
subsequent acquisitions and donations.

A library copy establishes availability, not reading.
Reviews, teaching, and identifiable uses provide evidence
of reception. Together, these sources place Bachelier's
later recognition against a background of early international
dissemination and selective engagement.

\section{Osborne and the problem of material accessibility}

The immediate stimulus for this investigation is a later
remark by M.~F.~M.~Osborne, a physicist at the U.S. Naval
Research Laboratory. In 1959, he published ``Brownian Motion
in the Stock Market'', an influential attempt to describe
stock-price fluctuations using methods inspired by statistical
physics \citep{Osborne1959}. He developed his approach
independently of Bachelier: the article contains no reference
to \emph{Théorie de la spéculation}.

After its publication, A.~G.~Laurent, a mathematician
working in statistics, then at Wayne State University
in Detroit, pointed out that related probabilistic
models of price movements had already been developed
and used in economic analysis \citep{Laurent1959}. Replying later
that year, Osborne explained that subsequent correspondence
had alerted him to earlier works unknown to him when
he wrote his article \citep{Osborne1959Reply}.
He singled out Bachelier's thesis in particularly strong
terms: ``the pioneer work on randomness in economic time
series, and yet most modern in viewpoint, is that of
Bachelier''. Osborne thus
recognized both Bachelier's priority and his strikingly
modern understanding of random economic fluctuations.

Osborne's discussion reveals close attention to several
central features of Bachelier's work: the probabilistic
treatment of price changes, the assumption of zero expected
gain in speculation---gains and losses balance on average---and
the Gaussian distribution of price changes, whose spread
increases with time. Because he considered the original
thesis ``rather inaccessible'', he devoted much of his reply
to explaining its contents. He also identified a location:
the rare-book room of the Library of Congress
\citep{Osborne1959Reply}.

The Library of Congress copy bears Bachelier's autograph
dedication to Robert S.~Woodward, a leading figure
in American mathematical and physical science and president
of the American Mathematical Society in 1899--1900.
A photograph of the cover of the separate Gauthier-Villars
edition of 1900, supplied to the author by Library of Congress
staff, clearly shows the dedication and signature.
Osborne did not mention the dedication; his reply focused
on Bachelier's mathematical argument and historical priority.

Four years later, Mandelbrot similarly described
Bachelier's original work as ``fairly inaccessible''
and supplied several bibliographical references,
including the thesis as published in the
\emph{Annales scientifiques de l'École normale supérieure}
\citep{Mandelbrot1963}.

Neither observation establishes how readily available
the thesis had been in 1900. They concern its accessibility
some six decades later: difficulty locating a work
in the late 1950s or early 1960s does not establish
difficulty obtaining it when it first circulated.

Osborne's experience thus connects two distinct questions.
Intellectual reception concerns who read Bachelier,
understood or cited his arguments, or incorporated his
results into later work. Material accessibility concerns
where copies travelled, which institutions received them,
how they were catalogued, and whether they remained
available. His rediscovery belongs to the history
of reception; his remark about inaccessibility prompts
an investigation of the circulation of copies.

Courtault and his coauthors subsequently emphasized
the thesis's publication both in a prestigious mathematical
journal and as a separate book, challenging the narrative
of complete neglect \citep{Courtault2000}.
Read alongside Osborne's remarks, this raises a concrete
question: how widely available was the thesis in the years
immediately following 1900? Answering it requires tracing
its channels of circulation and identifying institutions
where it could be consulted, independently of known
individual readers.

Later descriptions of the thesis as obscure, neglected, or 
forgotten rarely distinguish these two questions.

Several authors have attributed Bachelier's limited
reception to his subject. Robert W.~Dimand, a historian
of economic thought, and Colin Read, a professor of
economics and finance, emphasized the distance between
financial speculation and the usual concerns of
mathematicians and physicists
(\citealp{Dimand1993}; \citealp{Read2013}).
Philip Protter, a probabilist specializing in mathematical
finance, went further, characterizing Poincaré's attitude
as hostile and describing Bachelier's work as
``ignored and forgotten for some time''
\citep{Protter2001}. Jens Jackwerth, a finance scholar,
similarly wrote that ``his thesis was promptly forgotten
until the middle of the 20th century''
\citep{Jackwerth2004}.

The image of a neglected precursor also appeared in
\emph{The New Palgrave: Finance} (1989), edited by
John Eatwell, Murray Milgate, and Peter Newman.
This volume brought together articles on finance
selected from \emph{The New Palgrave: A Dictionary
of Economics} (1987), a major four-volume reference
work in economics. It thus presented Bachelier's
achievement within an authoritative reference collection
for economists.
As Dimand noted, its editors described finance as having
``its Gregor Mendel in Louis Bachelier''
\citep{Dimand1993}. The analogy combined
a foundational contribution with recognition delayed
until much later.

This narrative also appears in accounts for physicists
and mathematicians. Jeremy Bernstein, a theoretical
physicist and science writer, associated Bachelier's
rediscovery with economists in the 1950s
\citep{Bernstein2005} and described the thesis
as ``totally unknown to physicists''
\citep{Bernstein2006}.
Norman Biggs, a mathematician at the London School
of Economics, praised its originality while stating
that it ``remained unnoticed for many years''
\citep{Biggs2016}.
The physicists Peter Richmond, Jürgen Mimkes,
and Stefan Hutzler described the work
as ``published only in his PhD thesis'' and
``largely obscure for half a century''
\citep{RichmondMimkesHutzler2013},
leaving unmentioned its publication in a major
mathematical journal.

Davis and Etheridge distinguished mathematical recognition
from economic neglect: ``Among economists, Bachelier
was simply forgotten''
\citep{DavisEtheridge2006}.
More recently, Protter described the thesis's publication
as an event ``which should have ushered in a revolution
in thought for the twentieth century'', while observing
that it was largely ignored \citep{Protter2022}.
These assessments contrast its retrospectively recognized
importance with its limited early reception.

The following sections test these narratives against three
kinds of evidence: library records for material availability,
contemporary reviews for documented reading, and later
writings for citations and uses of Bachelier's work.

\section{The separate Gauthier-Villars thesis: publication, exchange, and international circulation}

Bachelier's \emph{Théorie de la spéculation} circulated in two
distinct physical forms: in volume~17 of the third series of
the \emph{Annales scientifiques de l'École normale supérieure}
in 1900, and as a separate seventy-page doctoral dissertation
issued by Gauthier-Villars
\citep{Bachelier1900,Bachelier1900Separate}.
These forms opened complementary routes of distribution.
The journal circulated through subscriptions and exchanges,
while the separate edition was distributed through university
dissertation exchanges and offered for sale through the book
trade. Personal presentations and subsequent transfers between
collections extended its circulation.

The separate edition's commercial availability is documented
in Gauthier-Villars's bulletin for the first quarter of 1900,
whose notice of the thesis is examined in the following
subsection.
Surviving copies must therefore be considered within a broader
history of commercial and institutional circulation:
institutions outside the Paris university exchange network
could purchase the thesis, even where individual acquisition
circumstances remain undocumented.

The separate edition also offered a potential source of
interest through its dedication to Henri Poincar\'e.
The brief tribute occupied a full page, prominently displaying
his name and the association with one of the period's leading
scientists \citep{Bachelier1900Separate}. The dedication was
absent from the journal version, which began directly with
the introduction \citep{Bachelier1900}.

A separate page at the end of the thesis edition stated the
subject assigned by the Faculty for Bachelier's second thesis:
the resistance of an unbounded viscous liquid, governed by
Navier's equations, to the small translational motions of an
immersed solid sphere. Alongside the study of speculation
stood a problem of mathematical physics involving viscous
resistance, closely related to the questions of particle
motion and diffusion that would soon become central to the
physical theory of Brownian motion. Bertrand Duplantier
highlighted this juxtaposition in 2005, while noting that
Bachelier's thesis on speculation made no explicit connection
between the two problems
\citep{Duplantier2005}.

\subsection{Commercial publication and promotion}

The \emph{Bulletin trimestriel des publications nouvelles de
la Librairie Gauthier-Villars} for the first quarter of 1900
presents Bachelier's thesis under ``New Publications''.
The notice and accompanying table of contents occupy almost
two pages, pp.~8--9. The entry describes a seventy-page volume,
dated 1900 and priced at five francs
\citep{GauthierVillars1900}.
A copy of this bulletin is preserved in the library
of the Museo Galileo--Istituto e Museo di Storia
della Scienza in Florence, an institution devoted
to the history of science. Library staff supplied
photographs of pages~8--9, making it possible to examine
the publisher's presentation in detail.
More than a brief catalogue entry, the presentation explains
the work's ambitions, methods, and results, emphasizing its
scientific and practical importance. It anticipates scepticism
about a mathematical theory of speculation, asking whether
one can ``formulate, so to speak, the determinism of chance''.
It explicitly identifies the intended audience:
``It is of equal interest to the mathematician and the banker''.
Speculators are also addressed as potential users.
The publisher thus presented the thesis to readers in both
mathematics and finance.

The notice follows the dissertation's argument in some detail:
stock exchange operations, the definition of the probabilities
under consideration, and the derivation of the probability law
by three methods. It singles out the ``law of the radiation
of probability'' as ``a true curiosity'', then discusses
applications to different transactions. Integral equations,
series expansions, and successive approximations appear
alongside comparisons between theoretical results and
observed prices.

The continuation addresses complex operations, the probability
that a price will reach a given level within a specified
interval, and the probable time at which it will do so.
The concluding assessment interprets agreement between
calculated and observed values as evidence that ``the market,
unbeknownst to itself, obeys a law that governs it:
the law of probability''. A detailed table of contents
identifies the questions and applications treated.

The unsigned text combines exposition with favourable
evaluation. Its authorship remains unidentified: it may have
been prepared by the publisher or a collaborator, or drawn
from material supplied by Bachelier. It nevertheless
establishes that Gauthier-Villars promoted the dissertation
through a substantial account of its scientific content,
allowing readers to learn about its principal questions
and claims without obtaining the book.

A documented purchase shows that the separate edition reached
a British scientific library through the book trade in 1900.
According to acquisition information supplied by the library
staff, the copy now held by the Science Museum Library
was purchased on 18 July from Williams and Norgate.
The library grew out of London's South Kensington Museum
and served students and the wider public seeking scientific
literature. Williams and Norgate were booksellers, publishers,
and importers specializing in scholarly works and foreign
scientific literature, with premises in London and Edinburgh.
The acquisition documents a commercial route from a Paris
publisher to an institutional readership abroad.

The separate edition was subsequently recorded in the
legal-deposit section of the \emph{Bibliographie de la France}
on 29 December 1900 \citep{BibliographieFrance1900}.
The notice identified the author, title, and publisher,
but gave no price. It offered librarians, booksellers,
and researchers another means of discovering the work.
This late-December registration did not mark the beginning
of circulation: both the publisher's bulletin and the London
purchase precede it.

\subsection{The university dissertation exchange system}

Commercial publication operated alongside an organized
system for distributing French doctoral dissertations.
The ministerial order of 21 July 1882 and Jules Ferry's
circular of 11 August 1882 concerning the
\emph{service des thèses} regulated the deposit
and distribution of dissertation copies
\citep{Ferry1882Arrete,Ferry1882Circulaire}.
Some copies were reserved for exchanges with foreign
universities; others circulated among French faculties
of the same discipline.

The foreign system was reciprocal: French universities
sent dissertations and academic publications abroad and
received corresponding material from participating
institutions. The Bibliothèque nationale's catalogue
of academic publications received through international
exchanges in 1900, published in 1901, contains entries
under thirty-eight foreign academic institutions
\citep[pp.~5--200]{BibliothequeNationale1901}.
These institutional headings provide a contemporary
record of participation in the exchange network.
Grouped by the political boundaries of the time,
the institutions were:

\begin{itemize}
    \item \textbf{Germany (21):} Berlin, Bonn, Breslau (now Wrocław, Poland), Erlangen, Freiburg im Breisgau, Giessen, Göttingen, Greifswald, Halle, Heidelberg, Jena, Kiel, Königsberg (now Kaliningrad, Russia), Leipzig, Marburg, Munich, Münster, Rostock, Strasbourg (now in France), Tübingen, and Würzburg.
    \item \textbf{Switzerland (6):} Basel, Bern, Fribourg, Geneva, Lausanne, and Zurich.
    \item \textbf{The Netherlands (3):} Groningen, Leiden, and Utrecht.
    \item \textbf{Belgium (2):} Ghent and Liège.
    \item \textbf{Sweden (2):} Lund and Uppsala.
    \item \textbf{Denmark (1):} Copenhagen.
    \item \textbf{The Russian Empire (1):} Jurjew, now Tartu in Estonia.
    \item \textbf{The United Kingdom (1):} Oxford.
    \item \textbf{The United States (1):} Johns Hopkins University in Baltimore.
\end{itemize}

Paris was thus connected to a broad European university network and, through Johns Hopkins, to the United States.

Within France, a separate mechanism provided for circulation
among faculties of the same order. For a Paris science
dissertation, this principally concerned other science
faculties and their university libraries. Information supplied
by the University of Bordeaux describes the practice as
distributing theses for preservation and dissemination to
French universities concerned with the corresponding
disciplines.

Bachelier's dissertation was therefore intended from the outset
to circulate beyond Paris, independently of individual demand
for a book on speculation. To determine whether this mechanism operated in practice,
the catalogues of the thirty-eight foreign institutions
were searched, and their libraries were contacted. The evidence includes dated acquisitions,
surviving copies in dissertation collections, and historical
catalogue records.

\subsection{Evidence of distribution through the exchange network}

Several institutions provide dated evidence of receipt.
At Göttingen, the thesis entered the library in 1901,
directly from the University of Paris through the university exchange system,
as the library established. It remains preserved with
other Paris dissertations.
Ghent likewise confirmed systematic receipt of dissertations
from French universities; its copy of Bachelier's thesis
was already present in 1900.

Strasbourg's copy is bound with dissertations from the Paris
Faculty of Sciences for 1899--1900. Strasbourg then belonged
to the German Empire, and the copy documents the thesis's
arrival in that academic library network.

Tübingen has confirmed that its copy survives.
The historical catalogue card describes the seventy-page Paris
edition of 1900 as
``Paris. Th. sc. math. 1899--1900. N° 1018.''

Erlangen confirmed contemporary receipt from its historical
card catalogue. Previously absent from the online catalogue,
the copy assigned to external storage received an online
record in response to this inquiry. The investigation thus
established historical receipt and made the holding
discoverable through the modern catalogue.

At Greifswald, another university on the Paris exchange list,
the library located a copy previously absent from the online
catalogue. Following its identification during this inquiry,
the library confirmed the creation of an online record.

Other German university collections retain copies or catalogue
records. Leipzig preserves the dissertation within its
mathematical holdings. Catalogued copies and records also
exist at Jena, Kiel, Freiburg im Breisgau, and LMU Munich.
Bonn, Halle, and Marburg also preserve copies of the
dissertation, documented by historical catalogue cards
that describe the Paris dissertation of 1900 under the
same designation as at Tübingen.

In Switzerland, Basel and Lausanne have confirmed surviving
copies. Lausanne's exemplar is now at the EPFL, in the
mathematical collection transferred from the University
of Lausanne in 2003. Basel's original remains in the library,
identified through the digitized historical dissertation
card catalogue rather than an individual modern online record.

The Royal Danish Library located an original copy in the
former dissertation collection of Copenhagen University
Library. Foreign theses there often received only minimal
cataloguing and were consequently absent from the modern
online catalogue. Around 1900, Copenhagen participated in
systematic dissertation exchanges with French institutions,
including the University of Paris. This inquiry prompted
the creation of an individual online record.

At Utrecht, Special Collections located an original copy
absent from the online catalogue and subsequently added it.

In Belgium, alongside Ghent's dated evidence, Liège's current
card catalogue identifies an available copy of the separate
dissertation of 1900. Johns Hopkins University, the only
American institution on the exchange list, holds a
catalogued exemplar.

Oxford provides another documented acquisition through
the Paris university exchange system. The Bodleian Libraries
preserve an original separate copy bearing an accession stamp
dated 6 April 1901. An accompanying University of Paris
exchange form, dated 30 November 1900, announces the dispatch
to Oxford of ``toutes les thèses et publications académiques''
for 1899--1900, totalling 658 items
(\emph{Bodleian Library Records}, c.~1044).

Lund preserves a copy of the separate Gauthier-Villars edition
of 1900, recorded in its old card catalogue.
A second copy is preserved in Mittag-Leffler's collection
at Lund's Centre for Mathematical Sciences.

Historical records also document copies no longer located.
At Jurjew, now Tartu, an old catalogue card establishes former
possession of the dissertation. Uppsala likewise retains
evidence of an unlocated copy. Surviving books thus represent
only part of the original distribution.

Subsequent management of these large collections complicates
the reconstruction of their contents. In 2024, Leiden University
Libraries transferred approximately 400,000 foreign dissertations
to the Internet Archive for digitization. Accumulated through
academic exchanges and largely uncatalogued, they occupied
more than 3.2 kilometres of shelving, creating substantial
storage pressures.\footnote{Leiden University Libraries,
annual report for 2024 and announcement concerning the
collection of foreign dissertations, September 2024.}
Library staff considered it plausible that Bachelier's
dissertation had belonged to this collection, but could
not confirm its presence.

Space constraints have also led to the disposal of older,
rarely consulted dissertations. Some libraries contacted
during this inquiry reported discarding large numbers
of such works together with their catalogue cards.
Present holdings and surviving records therefore reflect
both original circulation and subsequent decisions about
storage, description, and retention.

Staffing constraints also limited the investigation.
Some libraries reported that they currently lacked
the staff capacity to search their historical dissertation
collections and associated records. These cases remain
unresolved: the inability to undertake a search does not
establish the absence of Bachelier's thesis.

All fourteen institutions for which the inquiry obtained no confirmation of present or past possession appear in the catalogues of foreign dissertations received by the Bibliothèque nationale in both 1900 and 1901 \citep{BibliothequeNationale1901,BibliothequeNationale1902}. These records establish their active participation in the exchange network during the relevant period. Given the reciprocal organization of these exchanges and the documented receipt of Bachelier's thesis at other participating institutions, its reception at these fourteen establishments is probable, although it remains unconfirmed.

The French evidence complements these foreign holdings.
Bordeaux recorded its copy in 1900 and explicitly links it
to the organized distribution of theses among French
universities. Montpellier similarly identifies its exemplar
with systematic dissertation exchanges; it is preserved
among Paris doctoral theses in mathematical sciences.
The original Gauthier-Villars dissertation is also recorded
in the science libraries of Lyon, Poitiers, and Rennes,
extending its documented geographical presence in France.

Paris holdings include the Bibliothèque de la Sorbonne,
the Bibliothèque nationale de France, and the libraries
of the Institut Henri-Poincaré, the Institut de France,
the Muséum national d'histoire naturelle, the École normale
supérieure, and the École des mines de Paris.
Together, these holdings document preservation across
a broad range of French university and scientific libraries.

These findings show that the exchange list described an
effective mechanism of circulation, while revealing the
limits of searches confined to modern catalogues.
Basel, Lund, and T\"ubingen illustrate the continuing value
of historical cards. At Copenhagen, Greifswald, Erlangen, and Utrecht,
the inquiry led to the creation of individual online records. 

\subsection{Booksellers, personal presentations, and other collections}

The separate edition circulated beyond the university
exchange network. Its advertised price and the documented
London purchase establish its commercial availability.
Copies outside the exchange list belong to this ordinary
circulation of scientific books, alongside gifts and
transfers between collections. Their histories identify
several distinct routes.

Berkeley provides a clear example of bookseller-mediated
distribution. Its copy entered the University of California
library in November 1901 within a large, systematically
assembled collection of Paris mathematical dissertations
spanning much of the nineteenth century. Supplied from
Paris by the bookseller Hermann, the collection brought
the thesis to a university outside the official exchange
network.

The Library of Congress copy documents Bachelier's personal
initiative. Among the exemplars identified in this investigation,
it is the only known copy bearing his autograph presentation.
On the cover, he dedicated it to Robert S.~Woodward with
the words ``Hommage de l'auteur'' and his signature.
Woodward, a mathematician, geodesist, and physicist, had
presided over the American Mathematical Society in 1899--1900.
The dedication establishes that Bachelier addressed his work
directly to a prominent American scientist. A Library of
Congress accession stamp records its later entry in 1908.

Brown University's Sciences Library also catalogues a copy
of the separate Gauthier-Villars edition. Its provenance
and accession date remain unknown.

Another significant exemplar belongs to Gösta Mittag-Leffler's
personal collection of offprints, now at the Centre for
Mathematical Sciences in Lund. Founder of \emph{Acta Mathematica}
in 1882, Mittag-Leffler played a central role in international
mathematical communication, notably through close connections
with French and German mathematicians. The dissertation's
presence places it among the publications gathered by a
leading figure in these networks. Distinct from his book
library, the offprint collection was transferred from the
Mittag-Leffler Institute to Lund in 1992. Library staff
located the dissertation through its catalogue, compiled
in the 1990s.

A further connection with Mittag-Leffler's editorial activities
appears in the bibliographical supplement to volume~36 of
\emph{Acta Mathematica} (1913). On pages~2--3, the journal
listed Bachelier's \emph{Calcul des probabilités} (1912)
and summarized its contents \citep[pp.~2--3]{Acta1913},
explicitly including ``Théorie mathématique de la spéculation'',
``Étude des opérations de spéculation'', and
``Théorie de la spéculation''.

Other American copies have documented institutional histories
without an identified original supplier. Harvard's exemplar
entered the Harvard College Library in 1921, subsequently
joined Widener Library's economics holdings, and was transferred
to Baker Library in 1998. Its earlier provenance remains
unknown. At Cornell, marks on an original copy of the
separate edition of 1900 suggest acquisition in 1913. This would document continued
circulation more than a decade after publication, although
the marks do not identify the transaction.

In France, the Bibliothèque municipale de Lyon's copy came
from the Société de lecture de Lyon, a private learned
and cultural association with a library for its members.
This provenance places the thesis in a reading institution
outside the university system.

Italian copies reveal several further settings. The exemplar
at Ca' Foscari University's Biblioteca di Economia in Venice
belonged to the economist Alfonso de Pietri-Tonelli
(1883--1952), a professor there whose research included
mathematical economics and financial markets. Influenced
by Vilfredo Pareto, with whom he corresponded, he published
\emph{La speculazione di borsa}, a study of stock-market
speculation, in 1913. Bachelier's dissertation thus closely
matched his research interests. Its acquisition date and
source remain unknown, but another short publication by
Bachelier preserved with it bears the handwritten indication
``offert par l'auteur''. This documents an authorial
presentation within the collection, without establishing
that the dissertation itself was a gift.

The copy at the Institute for the History of the Contemporary
Age in Sesto San Giovanni, near Milan, came from the library
of the Collegio degli ingegneri e architetti di Milano
e provincia. This professional and learned association of
engineers and architects, with origins in 1563, was reorganized
in 1868 to promote technical knowledge and professional
practice. Its copy places the thesis in a setting devoted
to engineering and applied science.

At Sapienza University of Rome, the copy in the
``Enrico Barone'' Library belonged to the historical
collections of the institution preceding the modern
Faculty of Economics. Library research establishes its
presence by 1928, probably earlier.

Another Roman copy, in the historical library of the
Istituto Nazionale delle Assicurazioni (INA), came from
the mathematician Georg Bohlmann's private collection.
Its documented acquisition and former owner's scientific
interests warrant closer examination below.

These provenances document different routes: bookseller-mediated
acquisition, authorial presentation, and ownership by particular
institutions or individuals. Gauthier-Villars's offer for sale
provides the broader commercial context for acquisitions
beyond the university exchange system.

\subsection{From Bohlmann's private library to the INA in Rome}

The INA copy has a particularly informative provenance.
Founded in 1912, the INA was acquired in 2000 by
the Generali Group, which also took over its collections.
Generali, founded in Trieste in 1831, had developed
into a major international insurance group.
Its Historical Archive staff located the volume and
reconstructed its previous ownership and acquisition
through historical library records.

The volume entered the INA library on 20 February 1929
as part of a collection purchased from the Leipzig
antiquarian bookseller K.~F.~Koehlers. Comprising 90 rare
and out-of-print volumes and 1,100 pamphlets on pure
mathematics, probability, actuarial science, and statistics,
the collection came from the private library of the
mathematician Georg Bohlmann (1869--1928).

Bohlmann's career gives this provenance particular significance.
Supported by Felix Klein, he obtained his habilitation
at Göttingen in 1894 and taught actuarial mathematics
in the university's insurance seminar, established the
following year. In 1901 he received the title of extraordinary
professor, without a salaried professorship. In 1903 he
accepted a position with the Berlin branch of the New York
Mutual Life Insurance Company, while continuing his
mathematical research \citep[pp.~2--3]{Krengel2011}.

The acquisition documents the thesis's passage from a
specialist's personal library, through the antiquarian
book trade, to an Italian insurance institution.
Bohlmann's own acquisition date and circumstances remain
unknown, but the provenance places a surviving copy in
a mathematical collection directly relevant to Bachelier's
research.

A contemporary discussion supplies further context.
In November 1910, Bohlmann's and Bachelier's work appeared
in correspondence between the Russian probabilist
A.~A.~Markov and statistician A.~A.~Chuprov on the law
of large numbers for dependent variables. Responding
to Chuprov, Markov requested references to both authors'
articles. Bachelier's contribution was his 1906
``Théorie des probabilités continues''
\citep{Bachelier1906}, rather than
the 1900 thesis \citep[pp.~7--9]{Ondar1981}.
The exchange does not establish personal contact between
Bohlmann and Bachelier, but shows that their contributions
were considered together in contemporary discussions
of probabilistic dependence.

Bohlmann also contributed to the early axiomatization
of probability, in his work on life-insurance mathematics
published in 1900 and at the International Congress
of Mathematicians in Rome in 1908.
Already in 1900, Hilbert cited Bohlmann's
``Über Versicherungsmathematik'' in a footnote to his
sixth problem, which called for the axiomatic treatment
of physical sciences, particularly probability and
mechanics \citep[pp.~3--6]{Krengel2011}.

\subsection{The MIT copy and the later rediscovery}

The MIT case connects material circulation with later
intellectual rediscovery. Paul Samuelson (1915--2009),
an economist at the Massachusetts Institute of Technology
and a central figure in mathematical economics, received
the Nobel Memorial Prize in Economic Sciences in 1970.
His encounter with Bachelier brought the thesis to a
scholar whose work helped shape modern economic theory.

Samuelson knew of Bachelier before consulting the thesis.
In a recollection shared with Murad Taqqu on 14 August 2000,
he recalled the mathematician Stanislaw Ulam speaking
about Bachelier between 1937 and 1940, when both belonged
to Harvard's Society of Fellows. He also recalled
encountering Bachelier's name in William Feller's work
\citep[p.~27, including note~49]{Taqqu2001}.
This testimony identifies an earlier personal channel,
without specifying which works Ulam knew or what he
told Samuelson about them.

A later inquiry prompted Samuelson to seek Bachelier's
writing. He recalled that Leonard J.~Savage (1917--1971),
the American mathematician and statistician known for
\emph{The Foundations of Statistics} (1954), had sent
postcards to several economists, including himself,
asking whether they knew Bachelier and his popular
1914 book, \emph{Le jeu, la chance et le hasard}
\citep{Bachelier1914}. Searching the MIT mathematical
library for that book, Samuelson instead found the
1900 doctoral thesis. This consultation is associated
with the mid-1950s \citep[p.~vii]{Samuelson2006}.

The physical document remains unidentified. MIT's
separately catalogued item is a photocopy made around
1962 from an original at Princeton University Library.
It therefore cannot be the object Samuelson recalled
consulting several years earlier. MIT's current print
holdings also lack volume~17 of the third series of
the \emph{Annales scientifiques de l'École normale
supérieure}, containing the journal version.

Princeton's separate copy is bound as the first item
in the 1900 volume of \emph{Thèses Mathématiques}.
The digitized volume contains a bookplate identifying
Princeton's Mathematical Seminary. This was the original
used for MIT's photocopy.

Samuelson's account thus records access at MIT before
the surviving photocopy was produced. Princeton's copy
clarifies the later reproduction's history, but the
copy Samuelson consulted remains unknown.

\subsection{Donations of the thesis and subsequent publications}

The records of London's Institute of Actuaries document
another channel of circulation, reaching a professional
readership whose work connected probability with insurance
and financial valuation. The library additions published
in October 1906 list Bachelier's \emph{Théorie de la spéculation}
and \emph{Théorie mathématique du jeu}
\citep{Bachelier1901}, each marked
``From the Author'' \citep{Institute1906}.
This explicit identification of the donor establishes
Bachelier's personal role in presenting the thesis
to the Institute, six years after its publication.

Subsequent library lists record further donations
of his publications through 1913. Alongside the presentation
to Robert S.~Woodward, these repeated gifts document
Bachelier's sustained efforts to extend the international
circulation of his research and reach readers
in the actuarial profession.

A parallel donation is documented in Scotland.
Volume~XXIV of the \emph{Proceedings of the Edinburgh
Mathematical Society}, covering 1905--1906, lists
``Théorie de la Spéculation, Théorie Mathématique du Jeu,
par M.~L.~Bachelier'' among ``Presents to the Library''.
The introduction explicitly identifies the publications
as gifts. The record establishes receipt of the 1900
thesis and 1901 memoir by the end of that session
\citep{Edinburgh1906}: the same two works recorded
as donations to London's Institute of Actuaries in 1906.

The London Mathematical Society also recorded receipt
of \emph{Théorie de la spéculation} among the gifts
to its library between October 1905 and December 1906.
The thesis appears under the heading ``From the respective
Authors or Publishers'', which establishes receipt by gift
without identifying the donor more precisely
\citep{LMS1907}.
Subsequent library lists record further gifts of Bachelier's
publications.
\section{A second channel: circulation through the \emph{Annales}}

\emph{Théorie de la spéculation} also appeared in volume~17
of the \emph{Annales scientifiques de l'École normale
supérieure} \citep{Bachelier1900}. Founded in 1864 by
Louis Pasteur, the journal initially covered mathematics,
physics, chemistry, and the natural sciences. By around
1900, it had become devoted to mathematics, following
more than three decades of scientific publication associated
with one of France's leading institutions of higher learning.

The journal provided a second, independent channel of
diffusion: a library that never acquired or catalogued
the separate thesis could possess the complete text
through its subscription. Volume~17 circulated through
established networks of university and research libraries.
Its surviving copies offer another means of reconstructing
the geographical availability of Bachelier's work in the
years following publication.

\subsection{The early circulation of volume 17}

Libraries were contacted individually to establish whether
they hold an original copy of volume~17 of the third series,
published in 1900. Verification was necessary because
Kraus Reprint, based in Nendeln, Liechtenstein, reprinted
the volume in 1974. A catalogue entry dated 1900 could
therefore describe an original or a reproduction.
Library staff were asked to distinguish the editions
and, where possible, identify stamps, accession dates,
or other provenance evidence.

Bachelier's thesis appeared in the fascicle dated January
1900. The receipt dates discussed below document its
international circulation by April of that year, shortly
after his defence on 29 March. Receipt of an individual
fascicle could precede registration of the completed
annual volume.
The following overview groups holdings by country,
distinguishing confirmed originals from unresolved cases.

In Austria, original copies are held at the University
of Graz, the University of Innsbruck, the University
of Vienna, and the Vienna University of Technology.

In Belgium, original copies are preserved at the Catholic
University of Leuven, the Royal Library of Belgium
in Brussels, and the University of Liège.

In Britain, confirmed holdings include the British Library,
whose copy is stored at Boston Spa, Cambridge University
Library, the Bodleian Libraries at the University of Oxford,
the University of Edinburgh, the University of Glasgow,
and the University of Reading. Reading's volume formerly
belonged to the mathematician Eric Harold Neville, who
helped bring Ramanujan to Cambridge in 1914. The National Library of Scotland lists a run spanning 1866--1976, suggesting that its 1900 volume may be an original copy, although this remains unconfirmed.

In Denmark, original copies are held at the University
of Copenhagen's mathematics library and the University
of Southern Denmark in Odense.

In Germany, original copies have been confirmed at
Berlin State Library, Bielefeld University, TIB Hannover,
Göttingen SUB, Humboldt University Berlin, KIT Karlsruhe,
the Leibniz Institute for Astrophysics Potsdam,
Leipzig University, LMU Munich, and TU Munich;
at the university and state libraries in Bonn, Darmstadt,
Halle, and Münster; and at the universities of Cologne,
Erlangen--Nuremberg, Freiburg, Giessen, Greifswald,
Heidelberg, Jena, Rostock, and Tübingen.
Hamburg holds two reported copies: one at Hamburg
University's MIN-Forum Library, documented as a gift
received in 1922, and another confirmed at Hamburg
State and University Library.

In Italy, original copies are preserved at the National
Central Library of Rome, Sapienza University of Rome,
the University Library of Bologna, the University Library
of Genoa, the University of Genoa's science library,
the University of Florence, the University of Padua,
the University of Pavia, the Scuola Normale Superiore
in Pisa, the University of Turin, and the Library of
the Accademia Nazionale dei Lincei e Corsiniana in Rome.
At Sapienza, documentation from the booksellers Fratelli
Bocca records a subscription for 1900, providing evidence
of commercial circulation. In Naples, the Biblioteca Carlo Miranda and the Biblioteca della Società Nazionale di Scienze, Lettere e Arti hold extensive runs dating back to 1864, including uninterrupted sequences for 1896--1966 and 1884--1903, respectively. These holdings suggest that their 1900 volumes may be original copies, although neither library has confirmed this.

In the Netherlands, original copies are held at Delft
University of Technology, Leiden University, the University
of Amsterdam, and the University of Groningen.

In Spain, original copies are held at the Complutense
University of Madrid and the Library of Catalonia
in Barcelona.

In Sweden, original copies are preserved at Lund University,
Stockholm University, and Uppsala University.

In Switzerland, confirmed holdings include the Swiss
Federal Institute of Technology in Lausanne, the Swiss
Federal Institute of Technology in Zurich, the University
of Fribourg, and the University of Geneva.

In the United States, original copies are held at Brown
University, Bryn Mawr, Case Western Reserve, Columbia,
Cornell, Dartmouth, Harvard, Indiana University Bloomington,
Johns Hopkins, the Library of Congress, the New York Public
Library, Northwestern, Ohio State, Princeton, Stanford,
the University of California, Berkeley, the universities
of Chicago, Cincinnati, Illinois Urbana--Champaign, Iowa,
Kansas, Michigan, Minnesota, Missouri--Columbia, Penn State,
Pennsylvania, Pittsburgh, Wisconsin--Madison,
Wesleyan University, and Yale.

This list, based on correspondence with the institutions,
is not exhaustive. Nor does an original copy's present
location necessarily establish its presence there in 1900:
some volumes were acquired later or transferred.

Several records nevertheless document circulation within
months of publication. Bologna's chronological accession
register records a purchase on 2 April 1900 from Hoepli,
the Milanese bookseller and publisher, consistent with
receipt of the January fascicle. Oxford's volume bears
an accession stamp dated 17 April 1900, probably recording
receipt of that first issue. At Cornell, the preserved
January cover bears a receipt stamp dated 19 April;
Heidelberg received the fascicle on 21 April.
Freiburg's accession register also records arrival
in April 1900, and Uppsala confirms possession since
1900. These records document rapid circulation through
the journal's network, including arrival in North America
during the spring of 1900.

In France, the 1900 volume is recorded in university
libraries in Bordeaux, Grenoble, Lille, Lyon, Marseille,
Rennes, and Toulouse. Paris holdings include the École
normale supérieure, the Institut Henri-Poincaré, and
the mathematics research library at Jussieu.
Further holdings are recorded at the Bibliothèque
Jacques Hadamard in Orsay and the École polytechnique
in Palaiseau.

The extensive American holdings span the northeastern
states, the Midwest, and California, encompassing public
universities and privately endowed institutions.
This breadth should be viewed against the expansion
of American higher education in the second half of
the nineteenth century. Federal support, state investment,
and major private benefactions fostered institutions
across widely separated regions. The surviving copies
are consistent with the growth of this extensive network
collecting international scientific literature.

The numerous German holdings similarly reflect a dense
network of universities, technical institutions, and
research libraries, rooted in the university reforms
of the early nineteenth century.

\subsection{Encountering Bachelier through a mathematical journal}

The journal form also shaped how readers might encounter
Bachelier's work. While \emph{Théorie de la spéculation}
could suggest a specialized financial study, it appeared
in the \emph{Annales} among other contributions to a
major mathematical periodical.

The thesis occupied pages 21--86 of the January 1900
fascicle. The same annual volume
contained the French translation of David Hilbert's
\emph{Grundlagen der Geometrie}, first published in
German in 1899. Translated by Léonce Laugel as
\emph{Les principes fondamentaux de la géométrie},
it occupied pages 103--209, appearing across the March,
April, and May 1900 issues. Hilbert examined geometry's basic
assumptions and their logical relationships
\citep{Hilbert1900}.

The translation also contained additions supplied by
Hilbert. In a footnote on page 123, Laugel explained
that Hilbert had written a new note and an extended
addition to the conclusion specifically for the French
version. A note on page 204 identified the passage
extending to the table on page 208 as entirely new.
The translation thus offered supplementary material
absent from the German edition.

A reader consulting Hilbert's work, including these
additions, could therefore encounter Bachelier's study.
Published in different fascicles, the contributions
became part of the same annual volume preserved by
subscribing libraries. Bachelier's text reached
mathematical collections through ordinary journal
acquisitions, independently of specific interest
in financial questions.

Other contributions broadened the volume's potential
readership. Mechanics and mathematical physics occupied
a substantial place, with papers on rigid-body motion,
the elastic equilibrium of thin plates, the transverse
vibrations of elastic rods, and the equilibrium of an
ellipsoidal shell under internal pressure. Studies
in algebra, differential equations, and geometry
appeared alongside them.

The volume could thus attract mathematicians, mathematical
physicists, and engineers concerned with theoretical
mechanics, extending Bachelier's potential readership
beyond those deliberately seeking a study of speculation.

Approximately 160 copies can be traced through surviving holdings and historical records, indicating substantial international circulation. This evidence challenges the image of a lost or scarcely circulated thesis.

This early international availability supplies the context
for examining reception. Acquisition establishes access,
but does not itself demonstrate consultation or use.
Reviews, citations, teaching, and identifiable uses
of Bachelier's results document engagement, allowing
the relationship between circulation and intellectual
recognition to be examined through particular readers
and their responses.

\section{Early reviews and subsequent reception}

Bibliographical notices and reviews made the mathematical
content of Bachelier's thesis visible beyond institutions
holding copies. Readers could learn that
\emph{Théorie de la spéculation} offered a probabilistic
treatment of financial fluctuations. Later references
show that some incorporated its arguments and results
into their mathematical or economic writings.

\subsection{Contemporary notices and reviews, 1900--1902}

In 1900, Léon Autonne, of the Lyon Faculty of Sciences, reviewed Bachelier's thesis in the \emph{Revue générale des sciences pures et appliquées}, one of France's leading general scientific periodicals at the turn of the twentieth century. The journal combined original articles, scientific surveys, and bibliographical reviews, and published contributions by leading international scholars, including David Hilbert's famous address on mathematical problems given in 1900. Probably the first published review of \emph{Théorie de la spéculation}, the notice provides a contemporary assessment of the work \citep{Autonne1900}. Autonne briefly described its application to stock-exchange speculation and observed that appreciating the work required mathematical and financial expertise rarely found together. His notice thus identified the difficulty of assessing a contribution spanning these two fields.

A similar observation appeared in January 1901
in the Belgian \emph{Revue des questions scientifiques},
published by the Société scientifique de Bruxelles
\citep{RQS1901}. Unsigned, this notice described the
separate seventy-page volume issued by Gauthier-Villars
in 1900 and identified Bachelier's aim as uncovering
a probability law governing stock-market fluctuations.
Without assessing his conclusions, it noted that many
mathematicians would not be speculators, and many
speculators would not be mathematicians, presenting
this separation as a difficulty in evaluating and
using the results.

A more extensive discussion appeared in 1901 in
\emph{L'Enseignement mathématique}. Founded in 1899
by Henri Fehr in Geneva and Charles-Ange Laisant in
Paris, this international journal promoted exchanges
among mathematics teachers and connections with the
wider scientific community. Its committee included
prominent mathematicians from several countries.
A review could therefore introduce Bachelier's work
to an international readership engaged in mathematical
teaching and scholarship.

The reviewer, Louis Desaint, was a French mathematician
working in function theory. His dissertation,
\emph{Sur quelques points de la théorie des fonctions},
was presented at the Paris Faculty of Sciences in 1897
and published that year in the \emph{Annales scientifiques
de l'École normale supérieure}. His review thus records
the response of a researcher trained in mathematical
analysis.

Desaint devoted several pages to the thesis's financial
applications and mathematical arguments \citep{Desaint1901}.
He discussed the Gaussian probability law, the growth
of the characteristic size of fluctuations with the
square root of time, the connection between Bachelier's
``rayonnement de la probabilité'' and Fourier's heat
equation, and applications to stock-exchange contracts.

His assessment was explicitly favourable. Describing
the thesis as a ``beau travail scientifique'', he
congratulated Bachelier on his initiative and envisaged
an economic science increasingly informed by statistics
and mathematical analysis. His review shows a contemporary
reader recognizing central features of Bachelier's
construction and explaining them to a broader
mathematical audience.

The thesis also appeared in Eugène Estanave's 1901
\emph{Revue décennale des thèses}, a bibliography of
doctoral theses presented at the Paris Faculty of
Sciences between 1891 and 1900. Its entry identified
both the separate Gauthier-Villars edition and the
publication in the \emph{Annales}, and directed readers
to the reviews by Autonne and Desaint
\citep[p.~23]{Estanave1901}.
Without offering an assessment of its own, the
repertory provided a further means of locating
the thesis and contemporary discussions of it.

The work also entered German bibliographical channels.
In 1901, the \emph{Zeitschrift für Mathematik und Physik}
listed it as ``Theorie der Spekulation''
\citep{Wolffing1901}. This brief entry provided visibility,
though much less information than the extended reviews.

A fuller notice appeared on pages 241--242 of volume~31
of the \emph{Jahrbuch über die Fortschritte der Mathematik},
published in 1902 but covering publications from 1900.
It is signed ``Ot.'', identified in the collaborators'
list as Dr~Oster of Mannheim. Placed under combinatorics
and probability theory, it described the thesis as an ingenious memoir,
deserving attention from mathematicians and especially
readers concerned with finance and insurance.
Oster nevertheless stressed that stock-exchange terminology,
compounded by Bachelier's own definitions, hindered
explanation to non-specialists. His concise outline
covered mathematical expectation, probability laws,
probable deviations, options, and the probability of
reaching a specified price within a given interval.
It conveyed the investigation's breadth and recognized
its interest, without explaining the derivation connecting
the probabilistic model to Fourier's heat equation
\citep{Oster1902}.

These notices range from brief bibliographical entries
to detailed, favourable mathematical discussion.
They document recognition of the thesis's probabilistic
interest and difficulties in communicating it across
disciplinary boundaries.

\subsection{Subsequent references, 1908--1933}

Reception continued beyond these initial reviews.
Robert de Montessus de Ballore (1870--1937), a French
mathematician teaching at the Faculté libre des sciences
in Lille, worked on algebraic continued fractions,
earning recognition from the Académie des sciences.
His \emph{Leçons élémentaires sur le calcul des probabilités}
(1908) devoted a section to speculation and presented
results from Bachelier's thesis, explicitly referring
to a ``theorem of Bachelier'' on the speculator's zero
mathematical expectation \citep{Montessus1908}.
Incorporation into a probability textbook went beyond
announcing or evaluating the work.

Bernard Bru also identified the French actuary and
economist Alfred Barriol (1873--1959) as an expositor
of Bachelier's Gaussian model in \emph{Théorie et pratique
des opérations financières}, first published in 1908.
The accompanying note emphasizes its use by successive
generations of finance and insurance students
\citep[p.~16, note~19]{Taqqu2001}.
Yet Barriol did not acknowledge Bachelier: the calculations
appear without attribution in the first edition
\citep[pp.~324--327]{Barriol1908},
and \citet[p.~46]{BenElMechaiekhDimand2018} report
the same omission in the fourth
(1931, pp.~374--384).
His results could thus circulate without recognition
of their origin.

El Karoui, Parent, and Pradier situate Barriol's exposition
within actuarial teaching of mathematical economics.
They connect his presentation of Bachelier's methods
with his interest in Léon Walras (1834--1910), the
French economist who developed the mathematical theory
of general equilibrium. Barriol appears as an intermediary
bringing both authors' work to an actuarial readership
\citep{ElKarouiParentPradier2026}.

In 1910, Maurice Gherardt cited and discussed the thesis
in \emph{Le gain mathématique à la Bourse}, published
in Paris \citep{Gherardt1910}. A Paris Bourse speculator
and author of practical gambling manuals, he sought
to develop a profitable method of speculation.
His reference documents circulation among writers
addressing practising speculators.\footnote{Gherardt's
use of Bachelier's thesis had already been noted by
\citet{Jovanovic2000}. It was also discussed in Bernard
Bru's interview with Murad Taqqu
\citep[p.~16, including note~20]{Taqqu2001} and by
\citet[p.~46]{BenElMechaiekhDimand2018}.}

Another example is Alfonso de Pietri-Tonelli's
\emph{La speculazione di borsa}. An Italian economist,
Pietri-Tonelli taught commercial policy at the Royal
Higher Institute of Commercial Studies in Venice.
The revised and enlarged second edition (1921)
mentions Bachelier nine times
\citep{PietriTonelli1921}.\footnote{The Bocconi University
Library catalogue records an earlier edition published
in Rovigo by Tipografia sociale editrice in 1912--1913.
The references discussed here have been verified in
the second edition of 1921.}
The first mention cites \emph{Théorie de la spéculation}
(Paris, 1900) \citep[p.~4, n.~3]{PietriTonelli1921}.
Later references concern both the thesis and the
1912 \emph{Calcul des probabilités}, addressing
stock-exchange operations and mathematical expectation,
among other subjects. They show engagement beyond
an isolated bibliographical reference.

These references complement the material evidence:
Pietri-Tonelli's collection included the thesis now
at Ca' Foscari University of Venice. His published
references and surviving copy document both use
and possession of Bachelier's work.

In his 1931 paper on analytical methods in probability,
Andrei Kolmogorov identified Bachelier as, to his knowledge,
the first to study systematically schemes in which
transition probabilities vary continuously with time
\citep{Kolmogorov1931}. These probabilities describe
the likelihood of passing between states over a given
interval. He cited the 1900 thesis, Bachelier's 1910
paper on probabilities in several variables, and the
1912 \emph{Calcul des probabilités}
\citep{Bachelier1900,Bachelier1910,Bachelier1912}.

This acknowledgment has particular historical significance.
Kolmogorov became a leading twentieth-century mathematician
and a central figure in modern probability.
His 1931 paper helped establish the analytical theory
of Markov processes \citep{Kolmogorov1931}, while his
\emph{Grundbegriffe der Wahrscheinlichkeitsrechnung}
(1933) supplied a measure-theoretic axiomatic foundation
for probability \citep{Kolmogorov1933}.
His reference to Bachelier formed part of a foundational
contribution to the study of random processes
\citep[pp.~87--88]{DavisEtheridge2006}.

In 1933, Aleksandr Khinchine, a Russian mathematician
and major contributor to probability theory, also
acknowledged Bachelier. In a passage translated by
Taqqu, he credited him with pioneering the direct
study of distribution laws for continuous limiting
processes, while noting the inadequacy of his
mathematical tools \citep[p.~13, note~11]{Taqqu2001}.

These examples give precise meaning to selective
reception. Montessus incorporated thesis results
into a probability textbook; Barriol presented its
financial calculations to actuaries; Gherardt and
Pietri-Tonelli discussed it in writings on speculation.
Kolmogorov and Khinchine recognized its contribution
to probability laws evolving continuously with time.
Readers drew on different aspects according to their
questions and audiences. Barriol's exposition also
shows why explicit citations cannot fully measure
reception: results could be transmitted without
their author's name.

Against this background, the scarcity of references
among French mathematicians is striking. Beyond
the contemporary reviews by Autonne and Desaint,
Montessus's 1908 book remains the only example
identified here of a published discussion of the
thesis by a French mathematician before 1933.
This does not establish that it was unread in France,
but highlights the gap between availability and
visible reception in the mathematical literature examined.

The subject's cultural associations also deserve
consideration. In nineteenth-century France,
stock-exchange speculation provoked moral and political
controversy, particularly through its association
with gambling and rapid gains. Its legitimacy was
contested rather than universally rejected
\citep{Jovanovic2001}. Such associations may have
made speculation appear an unconventional subject
for mathematical research.

Reservations about Bachelier's manner of presenting
his results also help explain some mathematicians'
responses in the 1930s. In his conversation with
Murad Taqqu, Bernard Bru reports that Bohuslav Hostinský,
a Czech mathematician and theoretical physicist working
on probability, wrote to Maurice Fréchet, a French
mathematician known for his contributions to functional
analysis and probability theory, that Bachelier made
too many errors to be worth reading. Bru distinguishes this judgment
from the absence of rigorous proofs: Bachelier
often proceeded in the manner of a physicist,
indicating approaches and stating formulas rather
than supplying complete demonstrations.
In Bru's assessment, however, these approaches
and formulas were sound and frequently displayed
considerable mathematical originality
\citep{Taqqu2001}.
This testimony suggests that standards of exposition
and proof could hinder appreciation of Bachelier's
work, even when its results offered fruitful
directions for further research.

\subsection{Circulation, visibility, and mathematical use}

These findings support the reassessment by Courtault
et al. (2000) and Jovanovic (2002, 2012), who challenged
the image of Bachelier as a completely ignored precursor.
Jovanovic presented an unpublished manuscript by the
law professor Julien Laferrière, deposited in the Paris
Faculty of Law library in 1951, whose discussion of
the mathematical theory of speculation drew extensively
on Bachelier. It recorded thesis locations and shelfmarks
at the Bibliothèque nationale, the Paris Faculty of
Sciences, and the Bibliothèque Sainte-Geneviève, while
noting that the thesis and \emph{Calcul des probabilités}
had become difficult to obtain through booksellers
\citep{Jovanovic2002}.
This evidence connects intellectual reception with
the practical means of locating Bachelier's writings
before the rediscovery associated with Savage in the
1950s. The present inquiry complements this documented
readership by reconstructing the thesis's wider material
circulation. Mark H.~A.~Davis similarly described
Bachelier's work as a ``continuous thread'' through
twentieth-century stochastic analysis, contrasting
this mathematical reception with its much later uptake
by financial economists \citep[p.~9]{DavisBachelier}.

The early reviews further document contemporary
understanding and evaluation. Recognition of mathematical
interest need not have produced a broad readership.
The Belgian and German reviewers' difficulties connect
with Robert W.~Dimand's interpretation, which emphasized
the consequences of presenting a general mathematical
contribution through the specialized subject of financial
speculation \citep{Dimand1993}.

Davis and Etheridge contrast Bachelier's continued presence
in mathematical literature with his neglect among
economists \citep[pp.~2, 96]{DavisEtheridge2006}.
Gherardt's and Pietri-Tonelli's references qualify
this contrast: they establish that the thesis reached
writers on practical speculation and economic analysis,
without demonstrating a sustained research tradition
building on its financial applications.

\section{The significance of Bachelier's probabilistic construction}

The material and bibliographical evidence raises a further
question: what mathematical ideas circulated with Bachelier's
thesis? Although concerned with speculation on the Paris
Bourse, its conceptual reach extended beyond financial
valuation. Bachelier described price changes through
probability laws evolving over time, linking the composition
of successive fluctuations to a differential equation.

A central step was the transition from discrete trials
to continuous-time evolution. Classical problems of repeated
games typically indexed outcomes by the number of trials.
Continuous probability distributions were already well
established, but did not themselves describe random evolution
in continuous time. Bachelier made time a continuous parameter
in his description of fluctuations, a feature later emphasized
by Kolmogorov \citep{Kolmogorov1931}.

Bachelier assigned a probability distribution to the price
change over a given interval. He treated successive changes
as independent, with no preferred direction, and assumed
that their distribution depended on the interval's duration
rather than its starting date. The distribution over a
longer interval could therefore be obtained by composing
those over its constituent intervals.

In modernized notation, let $p(x,t)$ denote the probability
density of a price change $x$ over a duration $t$.
Bachelier's integral relation takes the form
\begin{equation}
p(z,t+s)
=
\int_{-\infty}^{+\infty}
p(x,t)\,p(z-x,s)\,dx,
\qquad t,s>0.
\label{eq:bachelier-composition}
\end{equation}
It expresses the independence of successive increments
and the dependence of their distributions on duration
alone \citep[pp.~35--36]{Bachelier1900}.

The probabilistic meaning is direct. A total change $z$
can result from any change $x$ during the first interval,
followed by a change $z-x$ during the second.
Independence allows the corresponding densities to be
multiplied; integration accounts for all intermediate
values $x$. Treating the interval as a whole must therefore
give the same distribution as dividing it into two parts,
whatever their durations.

The operation is the familiar convolution used to find
the density of a sum of independent random variables.
If $U$ and $V$ are independent random variables
with densities $f$ and $g$, then $U+V$ has density
\[
h(z)
=
\int_{-\infty}^{+\infty}
f(x)\,g(z-x)\,dx.
\]
In this usual application, $f$ and $g$ are given and
the formula determines $h$. It is not an integral
equation for an unknown family of densities, and
introduces neither time nor a probabilistic evolution.

Bachelier assigns convolution this further role.
The three densities in
equation~(\ref{eq:bachelier-composition}) belong to
the same unknown family, evaluated at $t$, $s$, and
$t+s$. The equation constrains how that family depends
on time. Bachelier explicitly calls it an equation
of condition that $p$ must satisfy
\citep[p.~36]{Bachelier1900}.
The requirement is consistency under every subdivision
of a time interval. In modern language, this equation is a special case
of the Chapman--Kolmogorov relation in which transition
probabilities depend only on the elapsed time and
the difference between the final and initial states
\citep{DavisEtheridge2006,Bernstein2005}.

Using this composition law to determine a continuous-time
probabilistic evolution was a pioneering contribution.
The novelty lay in the problem: finding a family of
probability laws describing change over arbitrary durations.
Kolmogorov's acknowledgment of Bachelier's systematic
study of transition probabilities depending continuously
on time helps identify its historical significance.

The distinction from statistical mechanics and kinetic
theory does not rest simply on a time-dependent density:
distributions in molecular physics could also evolve
continuously. It concerns the role of displacement laws
and their composition. In Bachelier's homogeneous model,
the same density describes a displacement over a given
duration from any initial price. Whereas Liouville's
equation describes transport of a phase-space distribution
by mechanical motion, and Boltzmann's equation describes
molecular transport and collisions, Bachelier imposes
consistency directly through assumptions about successive
random increments. His unknown family of displacement
laws must agree with its own composition over continuously
varying durations.

The connection to Fourier's equation can be understood by examining changes over a very short interval. The argument below is a modern, informal reconstruction: a rigorous derivation is mathematically delicate and requires more than the composition law alone.
Successive price changes are assumed to be independent, with distributions depending only on the duration of the interval. Their symmetry gives no preferred direction, and their mean square is assumed finite. Together, these properties imply that the mean square grows in proportion to elapsed time, without presupposing a Gaussian distribution.
Over an interval whose duration tends to zero, the displacement must become concentrated near zero, consistently with an initially vanishing price change. A stronger condition is also needed: changes exceeding any fixed size must contribute negligibly to the mean square per unit time as the interval shrinks. Finally, the density must be sufficiently regular to justify the expansions and limiting operations used below. Under these assumptions, the composition law leads to the diffusion equation.
Let $Y_s$ denote the price change during a short interval of duration $s$, with density $p(y,s)$. A bar denotes an average over possible changes, weighted by their probabilities. Under the assumptions stated above, the following limit exists and is finite:
\[
D=\frac{1}{2}\lim_{s\downarrow 0}
\frac{\overline{Y_s^2}}{s}
=
\frac{1}{2}\lim_{s\downarrow 0}
\frac{1}{s}
\int_{-\infty}^{+\infty}y^2p(y,s)\,dy.
\]
Since the mean increment is zero, its mean square is
also its variance. Thus $D$ measures half the rate at
which uncertainty accumulates with time. The expression
inside the limit is non-negative because $y^2$ and
$p(y,s)$ are non-negative and $s>0$. A larger $D$ means
more rapid spreading of the probability distribution.

To see how this coefficient enters the evolution
equation, write the composition law as
\[
p(x,t+s)
=
\int_{-\infty}^{+\infty}
p(x-y,t)\,p(y,s)\,dy.
\]

Although the integral extends over all possible changes, their probability becomes concentrated near zero as the interval shortens. The assumptions on larger changes and on the regularity of the density allow us to use the following local expansion inside the integral:
\[
p(x-y,t)
=
p(x,t)
-y\frac{\partial p}{\partial x}(x,t)
+\frac{y^2}{2}\frac{\partial^2p}{\partial x^2}(x,t)
+\cdots.
\]
On integration, the constant term gives $p(x,t)$,
because the density $p(y,s)$ integrates to one.
The linear term vanishes because $\overline{Y_s}=0$.
The quadratic term contributes
\[
\frac{\overline{Y_s^2}}{2}
\frac{\partial^2p}{\partial x^2}(x,t).
\]
Controlling the remainder is the delicate step in this argument. Under the assumptions stated above, its integrated contribution becomes negligible after division by the duration of the interval.
Since $\overline{Y_s^2}/s$ tends to $2D$, subtracting
$p(x,t)$, dividing by $s$, and letting $s$ tend to zero
gives
\[
\frac{\partial p}{\partial t}
=
D\,\frac{\partial^2p}{\partial x^2}.
\]

This is Fourier's heat equation, now governing
a probability density rather than temperature.
The analogy is instructive: just as heat spreads from
a concentrated source, probability spreads from an
initially concentrated distribution. The equation
thus describes a redistribution of probability,
without creating or destroying it.
This modern reconstruction connects the composition
law to the diffusion Bachelier also examined through
his ``radiation of probability''
\citep{Bachelier1900,Bernstein2005}.

Initially, the price change is exactly zero, so all
probability is concentrated at the origin. On the whole
real line, the normalized solution for $t>0$ is
\[
p(x,t)
=
\frac{1}{\sqrt{4\pi Dt}}
\exp\!\left(-\frac{x^2}{4Dt}\right).
\]
This Gaussian curve remains symmetric and centred
at zero, favouring neither positive nor negative changes.
Its total area remains one. Writing $X_t$ for the price
change over duration $t$, its mean and mean square are
\[
\overline{X_t}=0,
\qquad
\overline{X_t^2}=2Dt.
\]
Since the mean is zero, the variance is $\overline{X_t^2}$,
and the standard deviation is its square root:
\[
\sqrt{\overline{X_t^2}}=\sqrt{2Dt}.
\]
Variance grows linearly with time,
whereas the standard deviation, which measures the
characteristic size of fluctuations, grows as the square root of time. As the curve broadens,
its central height decreases as $t^{-1/2}$,
preserving the total probability.

The Gaussian solution satisfies the original
composition law. The composition law and the diffusion equation therefore
express the same evolution in two complementary ways:
one joins successive intervals, while the other describes
the change from one instant to the next.

Besides the integral approach, Bachelier obtained the
Gaussian law through a binomial approximation and through
his ``radiation of probability'', which leads to Fourier's
equation \citep[pp.~38--41, 45--49]{Bachelier1900}.
The integral formulation nevertheless has a wider scope
than this particular solution. The requirement that
successive intervals fit together consistently does
not itself require Gaussian changes. Different assumptions can produce
other probability laws while preserving the composition
principle.

The same idea can also accommodate a dependence on the
starting position or on time. In that case,
one must keep track of where and when each interval
begins, rather than only its duration and the size
of the change. This leads beyond the convolution formula
above to the broader Chapman--Kolmogorov framework.
What remains unchanged is the underlying requirement:
calculating a transition directly or through an
intermediate time must give the same probability.

Davis and Etheridge also locate the thesis's originality
in its treatment of trajectories as objects of probabilistic
reasoning \citep[pp.~6, 84, 90]{DavisEtheridge2006}.
Bachelier considered the probability of reaching
a specified price at some point during an interval.
This differs from asking where the price will be
at the interval's end. A price may reach a given level
and then fall back before the final observation.
The final distribution alone therefore does not answer
the question: the course of the price over time matters.

Bachelier did not construct a stochastic process on
a probability space in the modern sense. Kolmogorov
acknowledged his continuous-time schemes while explicitly
noting their lack of mathematical rigour
\citep{Kolmogorov1931}. Nevertheless, the thesis brought
together three related questions: how probabilities
combine across successive intervals, how their
distributions evolve with time, and what may happen
along an entire trajectory. Davis and Etheridge interpret
the subsequent rigorous development of these ideas
as a substantial part of the agenda of twentieth-century
probability theory.

\subsection{Smoluchowski, Chapman, and Kolmogorov:
different acknowledged antecedents}

The integral composition law connects three later
contributions: Marian Smoluchowski's study of Brownian
motion under external forces in 1913, Sydney Chapman's
study of Brownian displacements in 1928, and Andrei
Kolmogorov's analytical theory of stochastic processes
in 1931.
The quantities described differ across these settings:
price changes in Bachelier, particle positions and
displacements in Smoluchowski and Chapman, and
transitions between states of a system in Kolmogorov's
more general framework.

Smoluchowski, a Polish theoretical physicist and a pioneer
of the kinetic theory of matter, had already published
a study of Brownian motion in 1906. Using kinetic arguments
and a probabilistic analysis of successive molecular
collisions, he obtained the square-root-of-time law
for the characteristic magnitude of the displacement.
This derivation did not proceed through a transition
probability density evolving in continuous time,
an integral composition equation, or Fourier's diffusion
equation \citep{Smoluchowski1906}.

An explicit formulation of a composition law depending
separately on the initial and final positions appears
in Smoluchowski's 1913 study. The dependence on the
starting position arises from the presence of an external
force in his treatment of Brownian motion
\citep[pp.~423--425]{Smoluchowski1913}.

Smoluchowski names Einstein in the opening paragraph
of the 1913 paper, acknowledging his contribution to
the existing theory of Brownian motion, but gives no
specific attribution when deriving the integral
composition law. He does not mention Bachelier.
Einstein's earlier integral relation provides a plausible
source for this construction, although Smoluchowski
does not explicitly identify it as such.
The paper nevertheless documents an explicit
formulation of the transition composition law for
a spatially inhomogeneous case fifteen years before
Chapman's contribution.

Chapman, a British applied mathematician working notably
on the kinetic theory of gases, studied particle
displacements in a non-uniform fluid.
He opens by citing Einstein's Brownian-motion papers
of 1905 and 1906, presenting his study as an extension
from a uniform fluid to one whose temperature,
composition, or other properties affecting diffusion
vary with position \citep{Chapman1928}.
He does not cite Bachelier.

Both Smoluchowski's and Chapman's settings require
transition densities that can depend on the starting
position.
In Bachelier's model, $p(x,t)$ describes
the density of a price change $x$ over duration $t$.
Its form is independent of the initial price level,
giving the convolution in
equation~(\ref{eq:bachelier-composition}).
For the physical problems considered by Smoluchowski
and Chapman, the starting and ending positions must
instead be recorded separately.

To express the common mathematical structure,
consider a system whose state is described by a real
number. In modern notation, let $q_t(a,b)$ denote
the transition density from an initial state $a$
to a final state $b$ over duration $t$.
For a particle, these states are positions;
for a price model, they are price levels.
In a time-homogeneous setting, the density
satisfies
\begin{equation}
q_{t+s}(a,b)
=
\int_{-\infty}^{+\infty}
q_t(a,c)\,q_s(c,b)\,dc.
\label{eq:transition-composition}
\end{equation}
The integral accounts for every intermediate state $c$:
the system passes from $a$ to $c$, then from $c$ to $b$.
In Bachelier's model,
$q_t(a,b)=p(b-a,t)$, where $b-a$ is the price change.
Substitution therefore reduces this relation
to Bachelier's convolution.
For a particle in a homogeneous medium without
external forces, the same difference would represent
a spatial displacement.
The correspondence concerns the mathematical structure,
while the variables retain their distinct financial
or physical meanings.

Kolmogorov develops a general framework for transitions
between states, rather than restricting the variables
to prices or particle positions.
His formulation also allows transition probabilities
to depend separately on the initial and final times,
beyond the time-homogeneous notation used above.
In his introduction, he explicitly identifies Bachelier
as, to his knowledge, the first to investigate
systematically schemes with transition probabilities
varying continuously with time.
He cites \emph{Théorie de la spéculation} (1900),
\emph{Les probabilités à plusieurs variables} (1910),
and \emph{Calcul des probabilités} (1912)
\citep{Kolmogorov1931}.

The acknowledgment goes beyond this introductory
assessment.
Section~16 treats ``Bachelier's case'',
where, for given initial and final times, the transition
density depends on the difference between the final
and initial values of the state variable.
Kolmogorov attributes the corresponding forward
differential equation to Bachelier, distinguishing
its discovery from a rigorous proof.
Bachelier thus appears both as a historical antecedent
and as the author of specific cases incorporated
into the later theory.

Kolmogorov cites neither Einstein nor Chapman.
These papers therefore differ in their explicit
acknowledgments: Smoluchowski names Einstein when
introducing the existing Brownian-motion theory,
without specifically attributing the composition law;
Chapman presents his investigation as an extension
of Einstein's work; and Kolmogorov gives Bachelier
a prominent place.
These references document acknowledged antecedents,
but do not establish a complete chain of transmission
of the integral composition law.

\section{Later reception and renewed readership}

The early references establish Bachelier's visibility
in several settings. Later testimonies record how readers
encountered his work, their responses, and the place
they assigned it in probability and finance.
They reveal both the persistence of his ideas and
the surprise of finding, in a dissertation on speculation
from 1900, constructions associated with later developments.

Joseph Doob (1910--2004), the American mathematician
whose work helped establish modern probability,
particularly martingale theory, recalled his encounter
in an October 2003 communication to Davis and Etheridge.
When he began studying probability in 1934, references
in French texts led him to seek Bachelier's work.
Its ideas influenced his investigations of gambling
systems and subsequently martingale theory.
The recollection documents lasting influence, without
identifying either the French texts or the publications
by Bachelier he read \citep{DavisEtheridge2006}.

Recognition also entered mathematical terminology.
William Feller's influential probability treatise called
Brownian motion the ``Wiener--Bachelier process'', as Jeremy
Bernstein later noted \citep[p.~397, note~1]{Bernstein2005}.
This placed Bachelier within the teaching and exposition
of modern probability. Feller also credited him principally
with discovering the connection between random walks
and diffusion, emphasizing its importance for Kolmogorov's
theory of Markov processes, in a passage reproduced
by Taqqu \citep[p.~26, note~47]{Taqqu2001}.

Bachelier also occupies a prominent place in the opening
of Itô and McKean's \emph{Diffusion Processes and their
Sample Paths} (1965). Their preface credits him with
results describing Brownian motion, including its
connection with the heat equation and the probability
law of the maximum displacement over an interval.
They then discuss Einstein's physical derivation
of 1905, placing Bachelier's earlier contribution
within the historical development of the subject.
While observing that his ideas had not been appreciated
at the time, they distinguish his pioneering arguments
from the rigorous mathematical construction later
provided by Wiener \citep{ItoMcKean1965}.

This recognition also appears in university teaching.
In his 2017 lecture notes on probability and random
processes at KTH Royal Institute of Technology
in Stockholm, Timo Koski presents Bachelier as
the first to define the stochastic process now
called the Wiener process.
Introducing the subject, he highlights the thesis's
use of this process for option valuation and explicitly
recognizes Bachelier as a pioneer of both financial
mathematics and stochastic processes
\citep[p.~267]{Koski2017}.

Samuelson's encounter at MIT likewise illustrates
the intellectual impact of renewed access.
He called the thesis a ``gem'', valuing its tools
for describing uncertain changes over time and relating
their probabilities to equations for financial contracts.
He spread word of his discovery and encouraged wider
circulation, including publication in Paul Cootner's
1964 collection \citep[pp.~vii--viii]{Samuelson2006}.
In an interview reported by Davis and Etheridge,
he recalled what mattered most: ``it was the tools''.
The thesis offered ideas to reformulate and develop
in a new research setting.

Comparisons with Einstein's Brownian-motion research
became prominent in these retrospective assessments.
While emphasizing the soundness of Einstein's physics,
Samuelson argued that Bachelier had already obtained
overlapping mathematical results and explained how
Fourier's heat equation could describe probability
diffusion. This anticipation, together with temporal
composition, was central to his historical assessment
\citep[p.~viii]{Samuelson2006}.

Benoît Mandelbrot expressed similar astonishment
in \emph{Les objets fractals}. He emphasized the range
of later Brownian-motion mathematics anticipated
in a thesis on the Paris Bourse, whose wider significance
became apparent through subsequent developments
in probability \citep[pp.~170--174]{Mandelbrot1995}.

Mandelbrot reproduced extracts from Paul Lévy's letter
of 25 January 1964. A leading French probabilist
and major contributor to Brownian-motion theory,
Lévy discussed Bachelier's professional difficulties
and mathematical importance. While maintaining his
criticism of the initial formulation, he acknowledged
that Bachelier had identified important properties
before Einstein and Wiener, notably the diffusion
equation and the probability law of the maximum
over an interval. Lévy dated his reassessment to reading
Kolmogorov in 1931
(\citealp{Mandelbrot1995};
\citealp[p.~20]{Taqqu2001}).

In his autobiography, quoted separately by Taqqu,
Lévy stated that Einstein had rediscovered the connection
between probability density and the heat equation
in 1905, unaware of Bachelier's priority
\citep[p.~20, note~27]{Taqqu2001}.
These statements show how Lévy placed Bachelier
within Brownian-motion history, an assessment transmitted
by Mandelbrot and Taqqu to later readers.

The testimony also reached physicists.
Jeremy Bernstein acknowledged that Mandelbrot had
personally sent him a letter in which Lévy explained
his initial failure to appreciate Bachelier fully
\citep[p.~398, note~10]{Bernstein2005}.
Alongside the thesis, correspondence and recollections
thus shaped later assessments.

Mathematical appreciation did not extend equally
to the financial content. Courtault and his coauthors
report Lévy's annotation to his notebook summary:
``Too much on finance!'' \citep{Courtault2000}.
Recognition of the probabilistic results could coexist
with reservations about their financial setting.

Renewed readership can also be traced through Columbia
University's original volume~17 of the \emph{Annales
scientifiques de l'École normale supérieure}.
Containing the thesis, it had entered the collections
by 6 November 1901. Its later circulation record bears
Mandelbrot's name, Yorktown Heights, and a January 1961
date. The same volume connects early international
distribution with renewed research on financial
fluctuations. The entry documents borrowing, but does
not establish Mandelbrot's first encounter with Bachelier.

Mandelbrot's admiration accompanied criticism of the
financial model. Gaussian distributions assign very
small probabilities to changes far from the mean;
he argued that large movements occurred more frequently
in observed prices. He considered changes in logarithmic
prices, proposing distributions with heavier tails, allowing greater
probabilities of extreme changes.
He explicitly presented this as a generalization
of Bachelier's model
\citep[pp.~394--395]{Mandelbrot1963}.

He also called the Gaussian model both the prototype
and a limiting case of his theory, retaining its possible
relevance to economic phenomena with fewer extreme
observations \citep{Mandelbrot1963}.
Bachelier's construction was thus both an admired
achievement and a framework whose financial assumptions
could be questioned and extended.

Appreciation also reached financial practitioners.
Hélyette Geman emphasized the thesis's foundational
contribution to mathematical finance, recounting
the astonishment of traders shown its option-combination
payoff profiles from the discussion of ``opérations
complexes''. She also highlighted Robert C.~Merton's
acknowledgment of Bachelier. Merton, who received
the Nobel Memorial Prize in Economic Sciences in 1997, opened the
abstract of ``Theory of Rational Option Pricing''
(1973) by tracing option pricing to Bachelier's
1900 formula \citep[p.~41]{Geman1997}.

Merton returned to Bachelier's contribution in his
later discussion of mathematical models and financial
practice. He emphasized the richness
of Bachelier's thesis, highlighting mathematical
methods developed for speculation that anticipated
aspects of Einstein's later treatment of Brownian
motion \citep{Merton1995}.

Merton combined this recognition with the assertion
that Bachelier's work had remained unknown in
the financial literature for more than half a century.
His account thus illustrates how admiration for
the thesis could coexist with a narrative of prolonged
neglect. The early references and uses examined here,
particularly those of Barriol, Gherardt, and
Pietri-Tonelli, qualify that narrative by documenting
engagement with Bachelier's financial ideas before
their broader recognition in the mid-twentieth century.

Jeremy Bernstein offered a vivid assessment in
``Bachelier'', published in the \emph{American Journal
of Physics} in 2005. He presented the thesis as
anticipating the Brownian-motion theories of Einstein
and Smoluchowski \citep{Bernstein2005}.
Its Gaussian distribution would give a physicist
``an immediate sense of déjà vu''
\citep[p.~396]{Bernstein2005}.
His comparison examined probability composition,
square-root-of-time scaling, and Fourier's diffusion
equation, identifying specific mathematical correspondences
between the financial problem and later physical theories.

In ``Einstein and the existence of atoms'' the following
year, Bernstein argued that replacing price by position
revealed the mathematical problem Einstein later treated.
He identified the integral composition relation,
a random-walk approach comparable to Smoluchowski's,
and the Gaussian solution, even proposing the name
``Bachelier--Einstein equation''
\citep[pp.~869--870]{Bernstein2006}.
He had likewise described Einstein's integral relation
as a special case of the ``Chapman--Kolmogorov--Bachelier
equation'' \citep[p.~397]{Bernstein2005}.
These comparisons placed Bachelier within the history
of physical Brownian-motion theory.
His reference to Davis and Etheridge's 2006 translation
offered physicists another route to the thesis
\citep{Bernstein2006}.

Bertrand Duplantier likewise included Bachelier in his
historical account of Brownian motion. In his 2005 review,
he recalled Bachelier's use of a random-walk model in 1900
and situated his work in relation to the later mathematical
theory of Brownian motion \citep{Duplantier2005}.
This recognition extended beyond the biographical detail
of Bachelier's second thesis topic in hydrodynamics:
it gave his work on speculation a place in the history
of a central concept of physics. Duplantier's discussion,
like Jeremy Bernstein's, is noteworthy in a literature
where Bachelier's contribution has received uneven
attention. Whereas probabilists and mathematical
economists have explicitly acknowledged his mathematical
precedence, accounts by physicists and historians of
physics have more often centred on Einstein,
Smoluchowski, and Perrin, leaving Bachelier outside
the narrative. Duplantier and Bernstein thus stand out
for bringing his work into the historical perspective
of physics.

The economist and financial historian Peter L.~Bernstein
(1919--2009) devoted substantial attention to Bachelier
in his history of modern finance. He explained the
thesis's treatment of speculation as a fair game,
the speculator's zero mathematical expectation,
square-root-of-time fluctuations, and option valuation.
He also emphasized its anticipation of Einstein's
Brownian-motion research, presenting it as a major
achievement extending beyond its financial motivation
\citep{PeterBernstein1992}.

Bernstein contrasted this originality with limited
recognition, juxtaposing the distinguished career
it might have suggested with Bachelier's professional
difficulties. Drawing on Mandelbrot, he linked neglect
to the difficulty of classifying and retrieving findings
across disciplinary boundaries. His narrative culminated
in Savage's rediscovery in the 1950s, treating the thesis
as effectively lost \citep{PeterBernstein1992}.
The account illustrates the strength of the neglect
and rediscovery narrative, while highlighting a distinction
central here: limited recognition and difficult access
do not establish narrow circulation at publication.
Jeremy Bernstein later directed physicists to this
history of finance \citep{Bernstein2005}.

Tony Rothman, a physicist and science writer,
brought a similar comparison to a broader scientific
readership in ``Lost in Einstein's Shadow'',
published in \emph{American Scientist} in 2006.
He emphasized that Bachelier's financial
study produced an equation identical to Einstein's
later Brownian-motion equation
\citep[p.~113]{Rothman2006}.
Placing Bachelier among predecessors overshadowed
by Einstein's fame, Rothman also considered whether
Einstein knew his work, judging this unlikely because
of its supposed academic neglect.

Robert W.~Dimand had argued for Bachelier's place
in this history in \emph{The British Journal for
the History of Science} in 1993.
Responding to an account of Brownian motion, he stressed
the thesis's probabilistic construction and reproduced
Cootner's assessment of its continuous stochastic
processes, diffusion equation, first-passage problems,
and absorbing barriers \citep[p.~233]{Dimand1993}.
These appreciations situated Bachelier within probability
well beyond his financial application.

Dimand concluded that the same methods could describe
particle motion and financial price fluctuations.
He suggested that identifying the thesis too closely
with speculation had obscured its broader significance,
urging historians of Brownian motion to attend to it
\citep[p.~234]{Dimand1993}.

A public correction illustrates continuing reassessment.
In ``He Conceived the Mathematics of Roughness'',
his review of Benoît Mandelbrot's \emph{The Fractalist}
in \emph{The New York Review of Books} on 23 May 2013,
Jim Holt described Bachelier's model as borrowed from
the physics of a gas in equilibrium
\citep{Holt2013Review}.
Following readers' objections, he reconsidered this
in ``A Random Walk with Louis Bachelier'' on
24 October 2013. Acknowledging that he had understated
Bachelier's originality, Holt emphasized that financial
fluctuations had led him to Brownian-motion mathematics
five years before Einstein's physical investigation
\citep{Holt2013}.

Recent research has also shed new light on particular
arguments in Bachelier's thesis. Beiglböck, Pammer,
and Schachermayer (2022) revisit his idea of the
``radiation of probability'', in which probability
spreads between neighbouring price levels in a manner
analogous to the flow of heat.
The authors also draw a connection with modern optimal
transport, a field concerned with how one distribution
can be transformed into another by moving its mass \citep[sec.~1]{BeiglbockPammerSchachermayer2022}.
Their analysis illustrates how an argument developed
for financial speculation can continue to suggest
new connections within mathematics, extending the
thesis's reception beyond recognition of its pioneering
model of prices.

These episodes show different forms of renewed recognition:
Doob recalled an influence on his research; Samuelson
found tools opening new possibilities; Mandelbrot combined
admiration with criticism of the financial assumptions.
The Einstein comparisons by Samuelson, Lévy, both
Bernsteins, Rothman, and Holt were themselves part
of reception, expressing the importance they assigned
to mathematics developed through speculation.
Dimand likewise stressed its significance beyond
the original application.

Preserved copies, continuing references, personal
exchanges, and translations enabled these encounters.
Their history connects the thesis's early circulation
with the changing questions through which successive
readers appreciated it.

\subsection{English translations and wider circulation}

English translations extended the thesis's circulation
beyond readers of French and placed it in new intellectual
contexts. A.~James Boness's translation appeared in
Paul H.~Cootner's \emph{The Random Character of Stock
Market Prices} in 1964 \citep{Cootner1964}.
By including Bachelier alongside contemporary studies
of financial prices, the collection presented his work
as an antecedent to an active field of research.

Mark Davis and Alison Etheridge's 2006 edition offered
a different setting: a volume devoted to Bachelier's
thesis, combining a new translation with explanatory
annotations, a reproduction of the French original,
and an account of subsequent developments in probability
and financial economics \citep{DavisEtheridge2006}.
Whereas Cootner's collection situated the thesis among
studies of financial markets, Davis and Etheridge
provided a sustained examination of the text and its
place in the development of mathematical finance.

\section{Conclusion}

The evidence establishes early international circulation
through complementary channels. The separate dissertation
travelled through university exchanges, the book trade,
and personal gifts. Publication in the \emph{Annales
scientifiques de l'École normale supérieure} supplied
another route to European and American libraries.
The text reached institutions abroad soon after publication,
well before widespread recognition of its significance.

Bachelier actively contributed to this circulation.
The dedication to Woodward and repeated donations
to the Institute of Actuaries document his efforts
to reach readers abroad, including a professional
community concerned with probability and financial
valuation. The Institute's records trace this activity
from 1906 to 1913. The library lists of the London
Mathematical Society and the Edinburgh Mathematical
Society also document receipt of the thesis by gift
in the period 1905--1906.

Reviews and subsequent uses establish an identifiable
readership. Desaint discussed central mathematical
features; Montessus incorporated results into a probability
textbook; Barriol presented financial calculations
without attribution; Gherardt and Pietri-Tonelli discussed
the work in writings on speculation; and Kolmogorov
and Khinchine acknowledged its contribution to continuous-time
probability. These examples document selective reception
across intellectual and professional settings.
Together, these cases qualify
both the image of complete neglect and a sharp division
between mathematical recognition and economic neglect.
They establish dispersed engagement, without themselves
demonstrating a continuous research tradition.

The reconstruction also shows the value of combining
documentary sources. Present-day catalogues offer only
a partial view: surviving copies can remain difficult
to identify, while contemporary acquisition lists
record works no longer located. Volumes, provenance
marks, accession registers, historical catalogues,
and printed lists distinguish current holdings from
dated evidence of earlier availability.
This material and institutional inquiry extends previous
reassessments by specifying where, when, and through
which channels Bachelier reached potential readers.

Individual copies connect early acquisitions with later
encounters. Columbia's volume, accessioned in 1901,
records Mandelbrot's borrowing in 1961.
Samuelson's recollection of finding the thesis at MIT
likewise illustrates how institutional preservation
enabled later access. The translations of 1964
and 2006 widened access and supplied a new intellectual
context.

Delayed recognition must therefore be understood against
early international circulation, sustained efforts
by Bachelier to disseminate his research, and identifiable
but selective uses. The historical question is how
availability led to particular readings and applications,
and how these changed with mathematical questions
and readerships.

\section*{Acknowledgments}

The author is deeply grateful to the libraries and their staff whose generous and dedicated assistance made this inquiry possible. Their patient searches, careful examination of copies and historical records, and willingness to provide photographs and detailed information were essential to reconstructing the dissemination of Bachelier's thesis. Several went well beyond the initial request, locating previously uncatalogued copies and making them visible in online catalogues. This study owes much to their expertise, persistence, and remarkable spirit of collaboration.

\setlength{\bibsep}{3pt}

\end{document}